\documentclass[11pt]{article}

\newtheorem{thm}{Theorem}[section]

\newtheorem{lem}[thm]{Lemma}

\newtheorem{conj}[thm]{Conjecture}
\newtheorem{prob}[thm]{Problem}

\newcommand{\skipit}[1]{{}}
\newcommand{\prfend}{\hbox to7pt{\hfil}
\par\vskip-\baselineskip\hbox to\hsize
{\hfil\vbox {\hrule width6pt height6pt}}\vskip\baselineskip}

\newcommand{\myarrow}[2]{\hbox to #1pt{\hfil$\to$\hfil}{\hskip-#1pt{\raise 
10pt\hbox to#1pt{\hfil$\scriptscriptstyle #2$\hfil}}}}

\newcommand{\pr}[1]{\hbox{${\bf P}^{#1}$}}

\begin{document}

\title{The (unexpected) importance of knowing $\alpha$}

\author{Brian Harbourne\\
Department of Mathematics\\ 
University of Nebraska\\
Lincoln, NE 68588-0323 USA\\
email: bharbour@math.unl.edu}

\maketitle

\thanks{\noindent {\bf Acknowledgements}: The author thanks 
the sponsors of the Siena conference
for their support and for providing such a pleasant
environment for discussing mathematics, and the NSA for its support of his research.}

\begin{abstract}
We show that the determination, for all $r\ge0$,
of the Hilbert functions of all fat point subschemes of
\pr d with support at $r$ generic  points  is equivalent to determining,
for all $r\ge0$, the least degrees $\alpha$ such that
the Hilbert functions are positive (and hence to determining
the classes of all effective divisors on blow ups of
\pr d at $r$ generic points).
We also use this point of view for $d=2$ to show that the following 
conjecture is, surprisingly, equivalent to 
the standard conjecture for the Hilbert function of
fat points in the plane with generic support:
for any reduced irreducible curve $C$ on a blow up of
\pr2 at generic points, we conjecture that $C^2\ge g-1$, 
where $g$ is the arithmetic genus of $C$.
\end{abstract}

\section{Introduction}\label{intro}

For simplicity, we work over the complex numbers, ${\bf C}$. 
Fix an integer $d\ge 2$ (the case of $d=1$ being trivial).
Let ${\bf m}=(m_1,\ldots,m_r)$ be a finite sequence of
positive integers, let $p_1,\ldots,p_r$ be generic points of \pr d,
and consider the fat points ideal $I({\bf m},d)$
generated by all forms in ${\bf C}[x_0,\ldots,x_d]$ 
which vanish at each point $p_i$ to order at least $m_i$.

A significant open problem is: 

\begin{prob}\label{prb1}\rm
For each finite sequence
${\bf m}$ of positive integers, 
determine the Hilbert function of $I({\bf m},d)$; i.e., 
for each $t$, determine the dimension
$\hbox{dim}_{\bf C}I({\bf m},d)_t$ of the homogeneous
component of $I({\bf m},d)$ of degree $t$. 
\end{prob}

An apparently easier (but still open)
problem is:

\begin{prob}\label{prb2}\rm
For each finite sequence
${\bf m}$ of positive integers, determine
$\alpha({\bf m},d)$; i.e., 
the least degree $t$ such that
$\hbox{dim}_{\bf C}I({\bf m},d)_t>0$. 
\end{prob}

A related open problem is: 

\begin{prob}\label{prb3}\rm
For each integer $r>0$, determine the dimension
$h^0(X(r,d),\mathcal O_{X(r,d)}(D))$ of the complete
linear system of effective divisors
linearly equivalent to $D$, for every divisor
$D$ on the blow up
$X(r,d)$ of \pr d at $r$ generic points, $p_1,\ldots,p_r$.
\end{prob}

Another, seemingly easier, open problem is:

\begin{prob}\label{prb4}\rm
For each integer $r>0$, determine (inside
the divisor class group $\hbox{Cl}(X(r,d))$)
the subsemigroup $\hbox{EFF}(r,d)$ of 
classes of effective divisors on $X(r,d)$.
\end{prob}

The purpose of this note is to point out the not difficult
but perhaps unexpected and not widely appreciated fact 
that these problems are all equivalent.
Using this point of view, the
standard conjectural solution (versions of which have previously been given in
\cite{segre}, \cite{harb1}, \cite{Gthesis},
\cite{hirsch} and elsewhere)
to these open problems for $d=2$ can be reformulated
in a very concise way, as we show in Section \ref{theconj}.

\section{The Problems are Equivalent}

We first show that Problems \ref{prb1} 
and \ref{prb3} are equivalent.
Let ${\bf m}=(m_1,\ldots,m_r)$.
Then $\hbox{dim}_{\bf C}I({\bf m},d)_t$ is
equal to $h^0(X(r,d),\mathcal O_{X(r,d)}(D))$,
where $D=tE_0-m_1E_1-\cdots-m_rE_r$,
$E_0$ is the pullback to $X(r,d)$ from \pr d of the
class of a hyperplane, and for $i>0$, $E_i$ is the
class of the exceptional divisor corresponding 
to the blow up of the point $p_i$.
Thus a solution to Problem \ref{prb3}
implies a solution to Problem \ref{prb1}.
Conversely, given any divisor class $D$,
we have $D=tE_0-m_1E_1-\cdots-m_rE_r$
for some integers $t$ and $m_i$. 
By permuting the integers $m_i$, we 
obtain ${\bf m}'=(m'_1,\ldots,m'_r)$ with
$m'_1\ge\cdots\ge m'_r$. Let
$s$ be the greatest integer such that $m'_s>0$.
Note that $m_i<0$ implies that $-m_iE_i$ is in 
the base locus of $|D|$.
(To see this, take global sections of $0\to \mathcal O_{X(r,d)}(D-E_i)\to  
\mathcal O_{X(r,d)}(D)\to \mathcal O_{E_i}(m_i)\to 0$,
and induct.) We have 
$h^0(X(r,d),\mathcal O_{X(r,d)}(D))=
h^0(X(r,d),\mathcal O_{X(r,d)}(tE_0-m_1E_1-\cdots-m_rE_r))=
h^0(X(r,d),\mathcal O_{X(r,d)}(tE_0-(m_1)_+E_1-\cdots-(m_r)_+E_r))$,
where, for any integer $j$, we let $(j)_+$ denote
the maximum of $j$ and $0$. Since the points are generic,
we also have that 
$h^0(X(r,d),\mathcal O_{X(r,d)}(tE_0-(m_1)_+E_1-\cdots-(m_r)_+E_r))
=h^0(X(s,d),\mathcal O_{X(s,d)}(tE_0-m'_1E_1-\cdots-m'_sE_s))
=\hbox{dim}_{\bf C}I({\bf m}',d)_t$,
where ${\bf m}'=(m'_1,\ldots,m'_s)$.
Thus a solution to Problem \ref{prb1}
implies a solution to Problem \ref{prb3}.

Similarly, a solution to Problem \ref{prb4} implies
a solution to Problem \ref{prb2},
since, given ${\bf m}=(m_1,\ldots,m_r)$, 
$\alpha({\bf m},d)$ is the least $t$ such that
$h^0(X(r,d),\mathcal O_{X(r,d)}(tE_0-m_1E_1-\cdots-m_rE_r))>0$,
and hence such that $tE_0-m_1E_1-\cdots-m_rE_r$ is the class
of an effective divisor.
Conversely, $tE_0-m_1E_1-\cdots-m_rE_r$ is the class
of an effective divisor if and only if
$tE_0-(m_1)_+E_1-\cdots-(m_r)_+E_r$ is, hence
if and only if $t\ge\alpha({\bf m}',s)$, where
${\bf m}'$ is as above. Thus 
a solution to Problem \ref{prb2} implies
a solution to Problem \ref{prb4}.

It is now enough to check that a 
solution to Problem \ref{prb1} implies
a solution to Problem \ref{prb2}, and vice 
versa. Clearly, a solution to Problem \ref{prb1}
implies a solution to Problem \ref{prb2}.
Conversely, suppose we want to compute
$\hbox{dim}_{\bf C}I({\bf m},d)_t$.
Given an integer $i\ge0$, let ${\bf m}(i)$
denote the sequence $(m_1,\ldots,m_r,1,\ldots,1)$ with $i$
additional entries appended, each such additional entry equal to 1.
If $\alpha({\bf m},d)>t$, then clearly
$\hbox{dim}_{\bf C}I({\bf m},d)_t=0$.
Otherwise, let $j$ be the least integer such that
$\alpha({\bf m}(j),d)>t$. Then 
$\hbox{dim}_{\bf C}I({\bf m}(j),d)_t=0$,
but $\hbox{dim}_{\bf C}I({\bf m}(i),d)_t=
\hbox{dim}_{\bf C}I({\bf m},d)_t-i$ for $0\le i\le j$,
since imposing each single additional generic base point
to a nonempty linear system drops the dimension of the linear
system by exactly 1. Thus $\hbox{dim}_{\bf C}I({\bf m},d)_t=j$,
hence a solution to Problem \ref{prb2}
implies a solution to Problem \ref{prb1}.

\section{A New Formulation of the Standard Conjecture for \pr2}\label{theconj}

In this section we work on \pr2; i.e., we fix $d=2$.
We first recall the version of the standard conjectural
solution to Problem \ref{prb3} for \pr2, given
in \cite{harb1}.

\begin{conj}\label{conj1}
Let $X_r$ be the blow up of \pr2 at $r$ generic points.
For each $r\ge1$: if $C\subset X_r$ is a reduced, irreducible curve of 
negative self-intersection, then $C^2=-1$ and $C$ is smooth and rational;
moreover, if $D$ is an effective nef divisor on $X_r$, then $h^1(X_r,
\mathcal O_{X_r}(D))=0$.
\end{conj}

In this section we show this conjecture is equivalent to
the following:

\begin{conj}\label{conj2}
If $C\subset X_r$ is a reduced, irreducible curve
on the blow up $X_r$ of \pr2 at any $r\ge1$ generic points,
then $C^2\ge g-1$, where $g$ is the arithmetic genus of $C$.
\end{conj}

It will be useful here and later
to keep in mind that $h^2(X_r,\mathcal O_{X_r}(C))=0$
for any effective divisor $C$, or indeed for any divisor
$C$ (such as a nef divisor) such that $C\cdot E_0\ge 0$. (To see this,
recall that the canonical class on $X_r$ is 
$K_{X_r}=-3E_0+E_1+\cdots+E_r$. 
Since $E_0$ is the class of a line, it is nef, so
it follows from $E_0\cdot (K_{X_r}-C)=-3-C\cdot E_0<0$
that $K_{X_r}-C$ is not effective. Therefore
$h^2(X_r,\mathcal O_{X_r}(C))=h^0(X_r,\mathcal O_{X_r}(K_{X_r}-C))=0$.)
Thus, by Riemann-Roch, we have
$h^0(X_r, \mathcal O_{X_r}(C))=(C^2-C\cdot K_{X_r})/2+1+h^1(X_r, \mathcal O_{X_r}(C))$
for any divisor $C$ (such as an effective or nef divisor) with $E_0\cdot C\ge0$.

First we verify that Conjecture \ref{conj1}
implies Conjecture \ref{conj2}. 
Let $C$ be a reduced irreducible curve on $X_r$.
If $C^2\ge0$, then $C$ is effective and nef
so $h^1(X_r,\mathcal O_{X_r}(C))=0$ by Conjecture \ref{conj1}.
Thus $1\le h^0(X_r,\mathcal O_{X_r}(C))=(C^2-C\cdot K_{X_r})/2+1$,
so $C^2\ge C\cdot K_{X_r}$. Now, by the genus formula,
$2C^2\ge C^2+C\cdot K_{X_r}=2g-2$, so $C^2\ge g-1$.
If however $C^2<0$, then, by Conjecture \ref{conj1}, $C^2=-1$,
and $C$ is smooth and rational, so $g=0$. Thus $C^2=-1=g-1$,
so again $C^2\ge g-1$.

We now show, conversely, that Conjecture \ref{conj2}
implies Conjecture \ref{conj1}. If for any
reduced, irreducible curve $C$ we have $C^2\ge g-1$, then clearly
$C^2<0$ implies $C^2=-1$, $g=0$ and so $C$ is smooth and rational.
(It also follows that $C\cdot K_{X_r}=-1$ and that
$1=h^0(X_r, \mathcal O_{X_r}(C))\ge(C^2-C\cdot K_{X_r})/2+1=1$,
and hence $h^0(X_r, \mathcal O_{X_r}(C))=(C^2-C\cdot K_{X_r})/2+1$.)
So now it is enough to show that every effective nef divisor
$D$ on $X_r$ satisfies $h^1(X_r, \mathcal O_{X_r}(D))=0$, or,
what is by Riemann-Roch the same, that
$h^0(X_r, \mathcal O_{X_r}(D))=(D^2-D\cdot K_{X_r})/2+1$.
But this is a consequence of the following lemma.

\begin{lem}\label{fclem}
Let $D\subset X_r$ be an effective nef divisor.
Then Conjecture \ref{conj2} implies that one of the following holds:
\begin{enumerate}
\item $|D|=|lA|$ for some reduced irreducible smooth rational 
curve $A$ with $h^1(X_r, \mathcal O_{X_r}(lA))=0=A^2$ and $h^0(X_r, \mathcal O_{X_r}(lA))=l+1$;
\item $|D|=|lA|$ for some reduced irreducible divisor $A$ with $A^2=A\cdot K_{X_r}=0$,
$h^0(X_r, \mathcal O_{X_r}(lA))=1$ and $h^1(X_r, \mathcal O_{X_r}(lA))=0$; or
\item $|D|$ contains a reduced and irreducible member, and $h^1(X_r, \mathcal O_{X_r}(D))=0$.
\end{enumerate}
\end{lem}

{\em Proof.} First, consider the case that $D$ is reduced and irreducible.
By assumption, $D^2\ge g-1$, hence $2D^2\ge 2g-2=D^2+D\cdot K_{X_r}$, so
$D^2\ge D\cdot K_{X_r}$. Moreover, $D$ is nef, so $D^2\ge0$. 
If $h^0(X_r,\mathcal O_{X_r}(D))>1$,
then a general section of $|D-E|$ is still reduced and irreducible,
but $h^0(X_{r+1},\mathcal O_{X_{r+1}}(D-E))=h^0(X_r,\mathcal O_{X_{r}}(D))-1$,
where $E$ is the exceptional curve coming from
the blowing up $X_{r+1}\to X_r$ of an additional
generic point, and we identify $D$ with its pullback to $X_{r+1}$.
If $D^2=0$, then $(D-E)^2=-1$, hence as we saw above
$h^0(X_{r+1}, \mathcal O_{X_{r+1}}(D-E))=1$, and so
$2=h^0(X_r, \mathcal O_{X_r}(D))\ge(D^2-D\cdot K_{X_r})/2+1=2$,
so we have $h^0(X_r, \mathcal O_{X_r}(D))=(D^2-D\cdot K_{X_r})/2+1$,
as desired.
If $D^2>0$, then $(D-E)^2=D^2-1$,
$h^0(X_{r+1}, \mathcal O_{X_{r+1}}(D-E))=h^0(X_r, \mathcal O_{X_r}(D))-1$,
and $((D-E)^2-(D-E)\cdot K_{X_r})/2+1=(D^2-D\cdot K_{X_r})/2+1-1$.
Thus it is enough to show $h^0(X_{r+1}, \mathcal O_{X_{r+1}}(D-E))=
((D-E)^2-(D-E)\cdot K_{X_{r+1}})/2+1$. Continuing in this way,
we reduce to the case that $D^2=0$ (in which case we are,
as we have seen, done), or to the case that $D^2>0$ but 
$h^0(X_r, \mathcal O_{X_r}(D))=1$.
But in the latter case,
$1=h^0(X_r,\mathcal O_{X_r}(D))\ge(D^2-D\cdot K_{X_r})/2+1$,
hence $D\cdot K_{X_r}\ge D^2$, but $D^2\ge D\cdot K_{X_r}$,
so $D^2= D\cdot K_{X_r}$, and we again have 
$h^0(X_r,\mathcal O_{X_r}(D))=(D^2-D\cdot K_{X_r})/2+1$.

Thus $h^1(X_r, \mathcal O_{X_r}(D))=0$ if
$D$ is reduced and irreducible, so now assume that no member
of $|D|$ is reduced and irreducible. Then either:
\begin{description}
\item[(a)] $|D|$ has a fixed component but $D$ is not fixed;
\item[(b)] $|D|$ is fixed but $D$ is not reduced and irreducible; or
\item[(c)] $|D|$ is fixed component free, but its general section
is not irreducible, which by Bertini's theorem
means that $|D|$ is composed with a pencil.
\end{description}

Suppose that $|D|$ has a fixed component; let $N$ be a reduced
irreducible component of the fixed part of $|D|$, but assume $D\ne N$.
Choose a reduced irreducible component $A$ of the general
member of $|D-N|$. Then $D-(A+N)$ is 
effective, and we may assume either that $A^2\ge0$,
or that $A^2<0$ and hence $A$ is a fixed component
of $|D|$. 

First we show that $N^2\ge0$. Suppose $N^2<0$ (and hence
$N^2=-1=N\cdot K_{X_r}$). Since $D$ is nef,
$D$ must have a reduced irreducible component 
$A'$ meeting $N$ positively. As we saw above, $h^2$ vanishes, so 
$h^0(X_r,\mathcal O_{X_r}(A'+N))\ge((A'+N)^2-(A'+N)\cdot K_{X_r})/2+1
=h^0(X_r,\mathcal O_{X_r}(A'))+A'\cdot N>h^0(X_r,\mathcal O_{X_r}(A'))$, 
which contradicts $N$ being a fixed component.  Thus $0\le N^2$ and,
since $N$ is reduced and irreducible (so nef) and fixed, we have
$1=h^0(X_r,\mathcal O_{X_r}(N))=(N^2-N\cdot K_{X_r})/2+1$ and so
$N^2=N\cdot K_{X_r}$.

Since $N$ is nef, we see $N\cdot A\ge0$, but
$h^0(X_r,\mathcal O_{X_r}(A))=h^0(X_r,\mathcal O_{X_r}(A+N))
\ge((A+N)^2-(A+N)\cdot K_{X_r})/2+1
=h^0(X_r,\mathcal O_{X_r}(A))+A\cdot N$, so
it follows that $A\cdot N=0$. And now we see that
we cannot have $A^2<0$, since in that case $A$ is a fixed 
component, and the same argument we used for $N$ implies
that we would have $A^2\ge 0$. If $A^2>0$, then the
subspace orthogonal to $A$ must be negative definite
(by Sylvester's signature theorem and the Hodge index theorem;
see Remark V.1.9.1 of \cite{hart}),
which contradicts $N^2\ge0=A\cdot N$. Thus $A^2=0$.
The same argument with $A$ and $N$ switched
shows that $N^2=0$, and so also $-N\cdot K_{X_r}=0$. But 
$0<h^0(X_r,\mathcal O_{X_r}(A))=(A^2-A\cdot K_{X_r})/2+1$,
so $-A\cdot K_{X_r}\ge0$.

Now, since $N$ is nef,
it is {\it standard} \cite{harb2}; i.e., there is a birational
morphism $X_r\to \pr2$ and a corresponding exceptional configuration
$E'_0, E'_1,\ldots, E'_r$ such that
$N$ is a nonnegative integer linear combination
of the classes $H'_0,\ldots,H'_r$,
where $H'_0=E'_0$, $H'_1=E'_0-E'_1$, 
$H'_2=2E'_0-E'_1-E'_2$, and
$H'_i=3E'_0-E'_1-\cdots-E'_i$, for $i>2$.
Since $N\cdot H'_r=-N\cdot K_{X_r}=0$,
we see that $N\cdot H'_i\ge0$ for all $i$.
The only nontrivial nonnegative linear combinations $N$
of the $H'_i$ with $N^2=N\cdot H'_r=0$
are the nonnegative multiples of $H'_9$ (thus we see that
$r$ must be at least 9).
But  $h^0(X_r,\mathcal O_{X_r}(lH'_9))=1$ for all $l\ge0$, hence,
since $N$ is reduced and irreducible,
we have $N=H'_9$.

Now we show that $N=A$.
We have $0=N\cdot A=H'_9\cdot A\ge H'_r\cdot A=-K_{X_r}\cdot A\ge0$,
thus $E'_i\cdot A=0$ for all $i>9$, hence
$A$ is a linear combination
of $E'_0,\ldots,E'_9$, orthogonal to $H'_9$
with $A^2\ge0$. The only such classes
are the multiples of $H'_9$ itself (see, for example,
Lemma 2.2 of \cite{LH}). Since $A$ is reduced and irreducible,
we have $A=H'_9$, as before. Thus $D=lH'_9$ for some $l\ge2$,
and we have $h^0(X_r,\mathcal O_{X_r}(D))=1$
and $h^1(X_r,\mathcal O_{X_r}(D))=0$, giving
part (2) of the lemma. (This also shows that item (a)
above does not occur.)

So finally, suppose $D$ is fixed component free,
but does not have a reduced and irreducible
general member. Then it must be
composed with a pencil. Thus a general
member of $|D|$ is a sum $D_1+\cdots+D_l$
of reduced irreducible and linearly equivalent curves
(hence $|D|=|lD_1|$), with 
$2\le h^0(X_r,\mathcal O_{X_r}(D_1))=(D_1^2-D_1\cdot K_{X_r})/2+1$
(hence $2\le D_1^2-D_1\cdot K_{X_r}$),
and $h^0(X_r,\mathcal O_{X_r}(lD_1))\le l+1$.
Therefore, $l+1\ge h^0(X_r,\mathcal O_{X_r}(lD_1))\ge(l^2D_1^2-lD_1\cdot K_{X_r})/2+1$
(hence $2\ge lD_1^2-D_1\cdot K_{X_r}$ and so
$2\ge (l-1)D_1^2+D_1^2-D_1\cdot K_{X_r}\ge2+(l-1)D_1^2$,
which, since $D_1^2\ge0$, implies $D_1^2=0$ and so
$2=-D_1\cdot K_{X_r}$ and $g=0$). 
Now $l+1=h^0(X_r,\mathcal O_{X_r}(lD_1))\ge (l^2D_1^2-lD_1\cdot K_{X_r})/2+1
=l+1$, and part (1) follows.
\prfend

\end{document}